\documentclass[12pt,preprint]{aastex}

\setcitestyle{numbers}

\usepackage{graphicx}
\usepackage{epstopdf}

\makeatletter
\DeclareFontFamily{U}{tipa}{}
\DeclareFontShape{U}{tipa}{m}{n}{<->tipa10}{}
\newcommand{\arc@char}{{\usefont{U}{tipa}{m}{n}\symbol{62}}}%

\newcommand{\arc}[1]{\mathpalette\arc@arc{#1}}

\newcommand{\arc@arc}[2]{%
  \sbox0{$\m@th#1#2$}%
  \vbox{
    \hbox{\resizebox{\wd0}{\height}{\arc@char}}
    \nointerlineskip
    \box0
  }%
}
\makeatother

\begin{document}

\title{The Two Incenters of the Arbitrary Convex Quadrilateral}

\author{Nikolaos Dergiades and Dimitris M. Christodoulou}

\begin{abstract}
For an arbitrary convex quadrilateral $ABCD$ with area  ${\cal A}$ and perimeter $p$, we define two points $I_1, I_2$ on its Newton line that serve as incenters. These points are the centers of two circles with radii $r_1, r_2$ that are tangent to opposite sides of $ABCD$. We then prove that ${\cal A}=pr/2$, where $r$ is the harmonic mean of $r_1$ and $r_2$. We also investigate the special cases with $I_1\equiv I_2$ and/or $r_1=r_2$.
\end{abstract}

\section{Introduction}

We have recently shown [1] that many of the classical two-dimensional figures of Euclidean geometry satisfy the relation
\begin{equation}
{\cal A} = p r / 2 \ ,
\label{area} 
\end{equation}
where ${\cal A}$ is the area, $p$ is the perimeter, and $r$ is the inradius. For figures without an incircle (parallelograms, rectangles, trapezoids), the radius $r$ is the harmonic mean of the radii $r_1$ and $r_2$ of two internally tangent circles to opposite sides, that is
\begin{equation}
r = 2 r_1 r_2 / (r_1 + r_2) \ .
\label{harm} 
\end{equation}
Here we prove the same results for convex quadrilaterals with and without an incircle. These results were anticipated, but unexpectedly, the two tangent circles in the case without an incircle are not concentric, unlike in all the figures studied in previous work [1]. This is a surprising result because it implies that the arbitrary convex quadrilateral does not exhibit even this minor symmetry (a common incenter) in its properties, yet it satisfies equations~(\ref{area}) and~(\ref{harm}) by permitting two different incenters $I_1$ and $I_2$ on its Newton line for the radii $r_1$ and $r_2$, respectively. This unusual property of the convex quadrilateral prompted us to also investigate all the special cases with $I_1\equiv I_2$ and/or $r_1=r_2$.

\section{Arbitrary Convex Quadrilateral}

\begin{figure}
\begin{center}
    \leavevmode
    \includegraphics[trim=.1cm .1cm .1cm .1cm, clip, angle=0,width=15 cm]{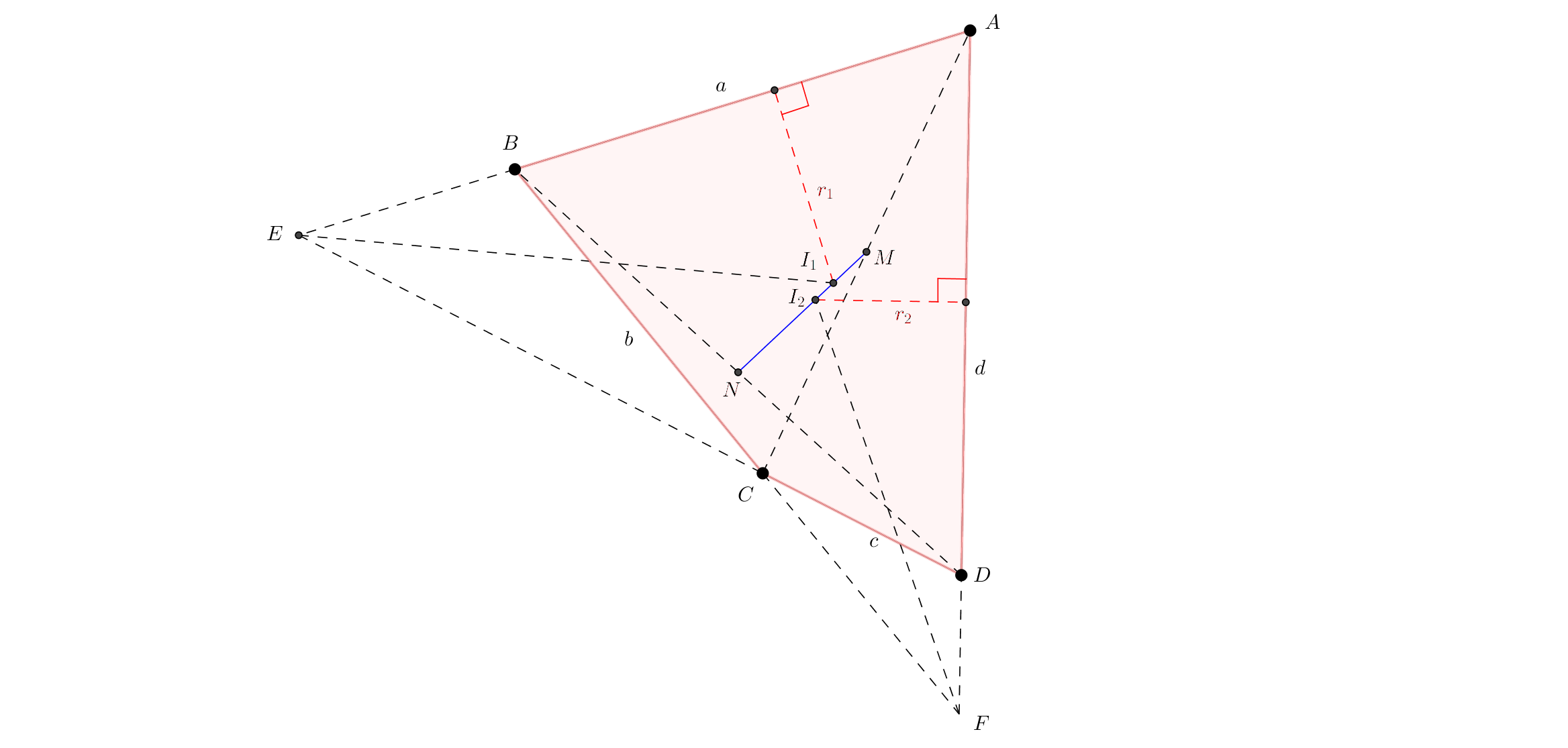}
\caption{Convex quadrilateral $ABCD$ with two interior incenters $I_1$ and $I_2$ on its Newton line $MN$.
\label{fig1}}
  \end{center}
  \vspace{-0.5cm}
\end{figure}

\begin{figure}
\begin{center}
    \leavevmode
    \includegraphics[trim=.1cm .1cm .1cm .1cm, clip, angle=0,width=15 cm]{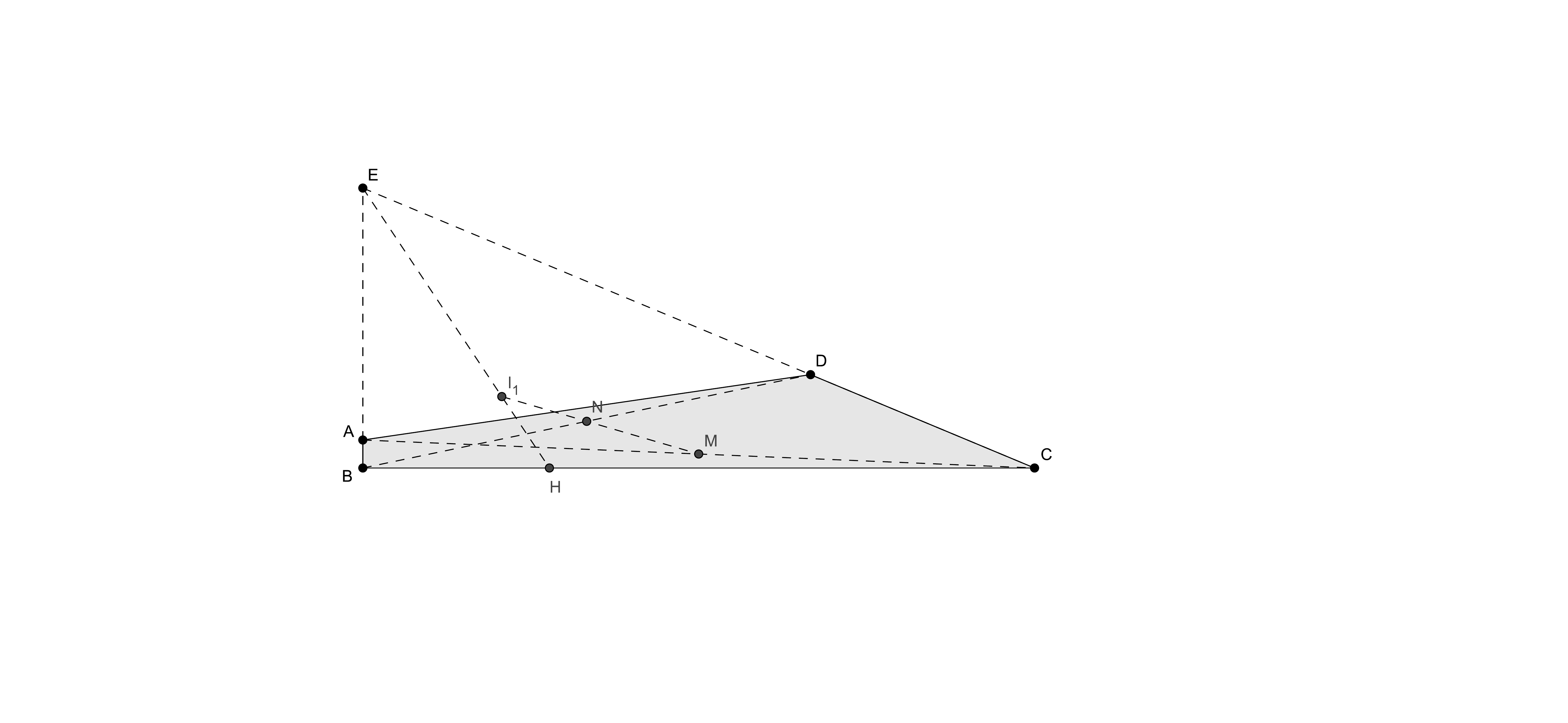}
\caption{Incenter $I_1$ lies outside of this quadrilateral $ABCD$.
\label{fig1_2}}
  \end{center}
  \vspace{-0.5cm}
\end{figure}

Consider an arbitrary convex quadrilateral $ABCD$ with Newton line $MN$ [2], where $M$ and $N$ are the midpoints of the diagonals $AC$ and $BD$ (Figure~\ref{fig1}).
Let the lengths of the sides of $ABCD$ be
$AB=a, BC=b, CD=c$, and $DA=d$. We extend sides $AB$ and $DC$ to a common point $E$. Similarly, we extend sides $AD$ and $BC$ to a common point $F$. We bisect $\angle E$ and $\angle F$. The angle bisectors $EI_1$ and $FI_2$ intersect the Newton line $MN$ at $I_1$ and $I_2$, respectively.

\noindent
{\bf Definition.} We define $I_1$ as the incenter of $ABCD$ that is equidistant from sides $AB$ and $CD$ at a distance of $r_1$. We also define $I_2$ as the incenter of $ABCD$ that is equidistant from sides $BC$ and $DA$ at a distance of $r_2$.

\noindent
{\bf Remark 1.} Points $I_1, I_2$ are usually interior to $ABCD$, but one of them can also be outside of $ABCD$ (as in Figure~\ref{fig1_2}).

\noindent
{\bf Lemma 1 (Based on L\'eon Anne's Theorem \#555 [3]).} Let $ABCD$ be a quadrilateral with $M, N$ the midpoints of its diagonals $AC, BD$, respectively. A point $O$ satisfies the equality of areas
\begin{equation}
(OAB) + (OCD) = (OBC) + (ODA) \, ,
\label{Oa1} 
\end{equation}
if and only if $O$ lies on the Newton line $MN$.

\noindent
{\bf Proof.} 
Using the cross products of the vectors of the sides of $ABCD$, equation~(\ref{Oa1}) implies that
\begin{equation}
 \begin{array}{ll} 
(OAB) - (OBC) + (OCD) - (ODA) = 0  & \iff\nonumber    \\
\overrightarrow{OA}\times\overrightarrow{OB} + \overrightarrow{OC}\times\overrightarrow{OB} +
\overrightarrow{OC}\times\overrightarrow{OD} + \overrightarrow{OA}\times\overrightarrow{OD} = \overrightarrow{0} & \iff\nonumber \\
(\overrightarrow{OA} + \overrightarrow{OC})\times\overrightarrow{OB} +
(\overrightarrow{OC} + \overrightarrow{OA})\times\overrightarrow{OD} = \overrightarrow{0} & \iff\nonumber \\
(\overrightarrow{OA} + \overrightarrow{OC})\times (\overrightarrow{OB} + \overrightarrow{OD}) = \overrightarrow{0} & \iff\nonumber \\
2\overrightarrow{OM}\times 2\overrightarrow{ON} = \overrightarrow{0} , \nonumber
         \end{array} 
\label{lemma1}
\end{equation}
therefore point $O$ lies on the line $MN$ (see also [4]).\hfill\sq

Since for signed areas it holds that $(OAB) + (OBC) + (OCD) + (ODA) = (ABCD)$, we readily prove the following theorem:

\noindent
{\bf Theorem 2 (Arbitrary Convex Quadrilateral).}  
The area of $ABCD$ is given by equation~(\ref{area}), where the radius $r$ is given by equation~(\ref{harm}) and the two internally tangent circles to opposite sides $\odot I_1$ and $\odot I_2$ are centered on two different points on the Newton line $MN$. 

\noindent
{\bf Proof.} 
Since $I_1$ lies on the Newton line, we find for the area ${\cal A}$ of $ABCD$ that
\begin{equation}
{\cal A}/2 = (I_1AB) + (I_1CD) = (a + c) r_1 / 2 \, ,
\label{la2} 
\end{equation}
or
\begin{equation}
a + c = {\cal A}/r_1 \, .
\label{la3} 
\end{equation}
Similarly, we find for the incenter $I_2$ that
\begin{equation}
b + d = {\cal A}/r_2 \, .
\label{la4} 
\end{equation}
Adding the last two equations and using the definition of perimeter $p=a+b+c+d$, we find that
\begin{equation}
p = {\cal A}\left(\frac{1}{r_1} + \frac{1}{r_2}\right) = {\cal A}~\frac{r_1+r_2}{r_1 r_2} \, ,
\label{la5} 
\end{equation}
which is equation~(\ref{area}) with $r$ given by equation~(\ref{harm}).
\hfill\sq

\begin{figure}
\begin{center}
    \leavevmode
    \includegraphics[trim=.1cm .1cm .1cm .1cm, clip, angle=0,width=15 cm]{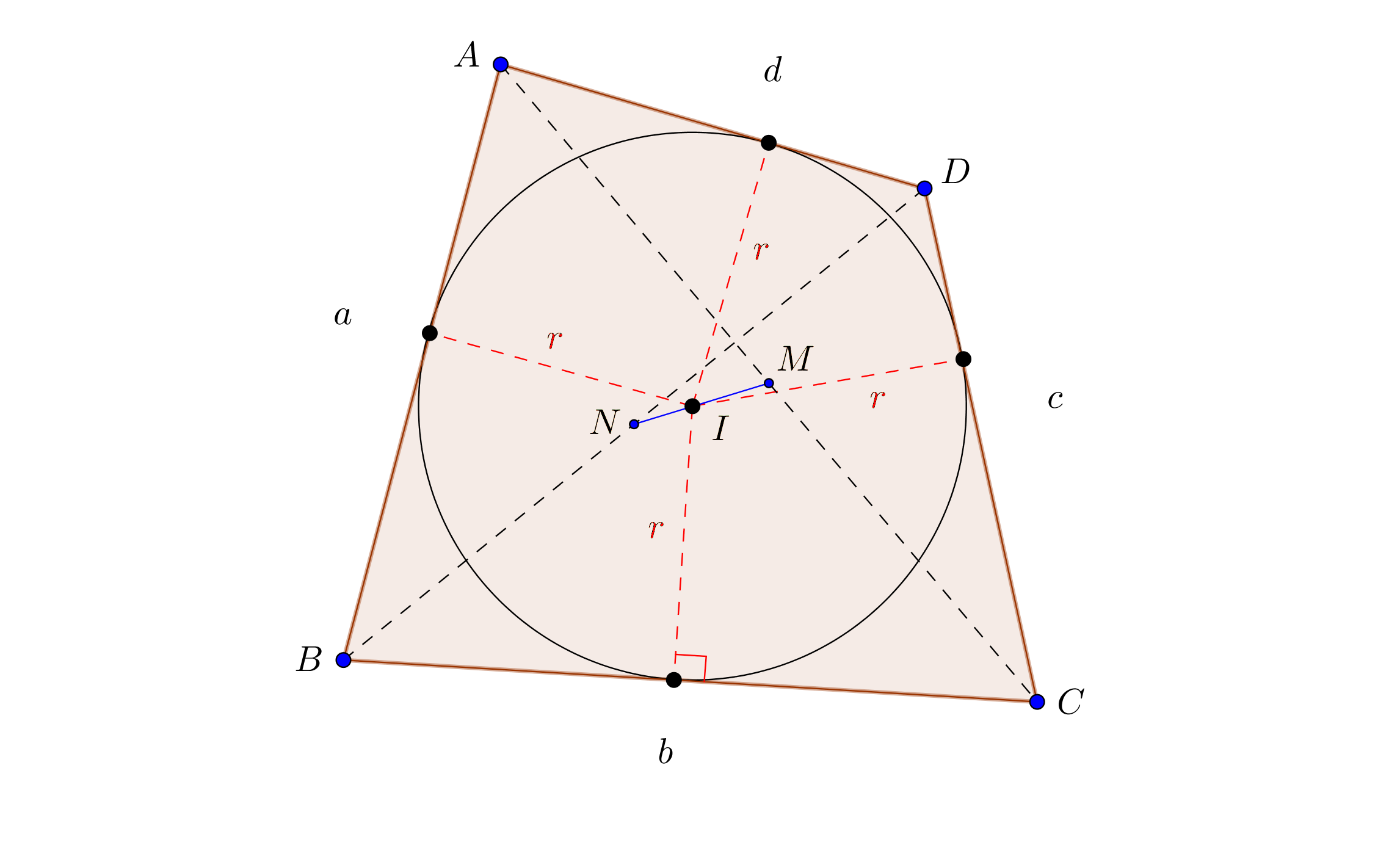}
\caption{Tangential quadrilateral $ABCD$ with a single incenter $I$ on its Newton line $MN$.
\label{fig2}}
  \end{center}
  \vspace{-0.5cm}
\end{figure}

\section{Tangential Quadrilateral}

In the special case of a tangential quadrilateral, the two incenters $I_1$ and $I_2$ coincide with point $I$ and, obviously, $r_1=r_2$. Using Lemma 1, we prove the following theorem:

\noindent
{\bf Theorem 3 (Based on Newton's Theorem \# 556 [3]).} If a quadrilateral $ABCD$ is tangential with incenter $I$ and inradius $r$, then $I$ lies on the Newton line $MN$ (as in Figure~\ref{fig2}), $I_1\equiv I_2\equiv I$, $r_1=r_2=r$, and the area of the figure is given by ${\cal A}=pr/2$.

\noindent
{\bf Proof.} If $ABCD$ has an incircle $\odot I$ of radius $r$, then by the tangency of its sides 
\begin{equation}
AB + CD = BC + DA \, .
\label{in1} 
\end{equation}
Multiplying by $r/2$ across equation~(\ref{in1}), we find that 
\begin{equation}
(IAB) + (ICD) = (IBC) + (IDA) \, ,
\label{in2} 
\end{equation}
where again parentheses denote the areas of the corresponding triangles. Therefore, by Lemma 1, the incenter $I$ lies on the Newton line $MN$ of $ABCD$ and the area of the figure is 
\begin{equation}
{\cal A} = (a + c) r / 2 + (b + d) r / 2 = p r / 2\, .
\label{in3} 
\end{equation}
Since $I$ is the point of intersection of the bisectors of $\angle E$ and $\angle F$ (Figure~\ref{fig2}), it follows that $I_1\equiv I_2\equiv I$, $r_1=r_2=r$.
\hfill\sq

\noindent
{\bf Remark 2.} The proof of the converse of Theorem 3 is trivial: If $I_1\equiv I_2\equiv I$ and $r_1=r_2=r$, then equation~(\ref{in3}) implies equation~(\ref{in2}) which in turn implies equation~(\ref{in1}).

\noindent
{\bf Theorem 4.} In quadrilateral $ABCD$, if the incircle $I_1(r_1)$ is tangent to a third side, then $ABCD$ is tangential.

\begin{figure}
\begin{center}
    \leavevmode
    \includegraphics[trim=.1cm .1cm .1cm .1cm, clip, angle=0,width=15 cm]{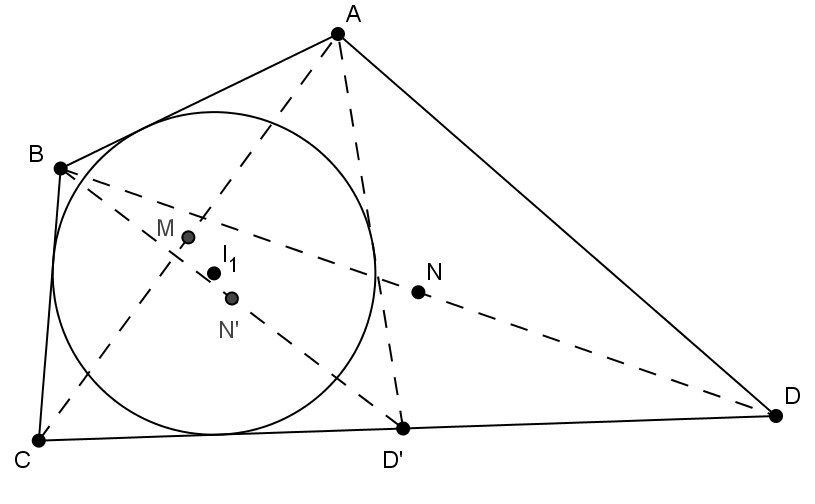}
\caption{Quadrilateral with an incircle $I_1$ tangent to three of its sides.
\label{newfig3}}
  \end{center}
  \vspace{-0.5cm}
\end{figure}

\noindent
{\bf Proof.} Let the incircle be tangent to side $BC$ in addition to sides $AB$ and $CD$ (as in Figure~\ref{newfig3}). If it were not tangent to side $AD$ as well, then we would draw another segment $AD^\prime$ tangent to the incircle and hence the incenter $I_1$ that lies on the Newton line $MN$ would also lie on the ``Newton line'' $MN^\prime$ of the tangential quadrilateral $ABCD^\prime$. Since $AB, CD$ are not parallel, the two Newton lines do not coincide, hence $I_1\equiv M$. In a similar fashion, if we draw a tangent line to the incircle from vertex $D$, we find that $I_1\equiv N$. But $M, N$ cannot coincide because $ABCD$ is not a parallelogram, thus the two equivalences of $I_1$ are impossible, in which case the incircle $I_1(r_1)$ must necessarily be tangent to the fourth side $AD$, making $ABCD$ a tangential quadrilateral. The same holds true for the incircle $I_2(r_2)$ when it is tangent to three sides of $ABCD$.
\hfill\sq

\section{Cyclic Quadrilateral}

\begin{figure}
\begin{center}
    \leavevmode
    \includegraphics[trim=.1cm .1cm .1cm .1cm, clip, angle=0,width=15 cm]{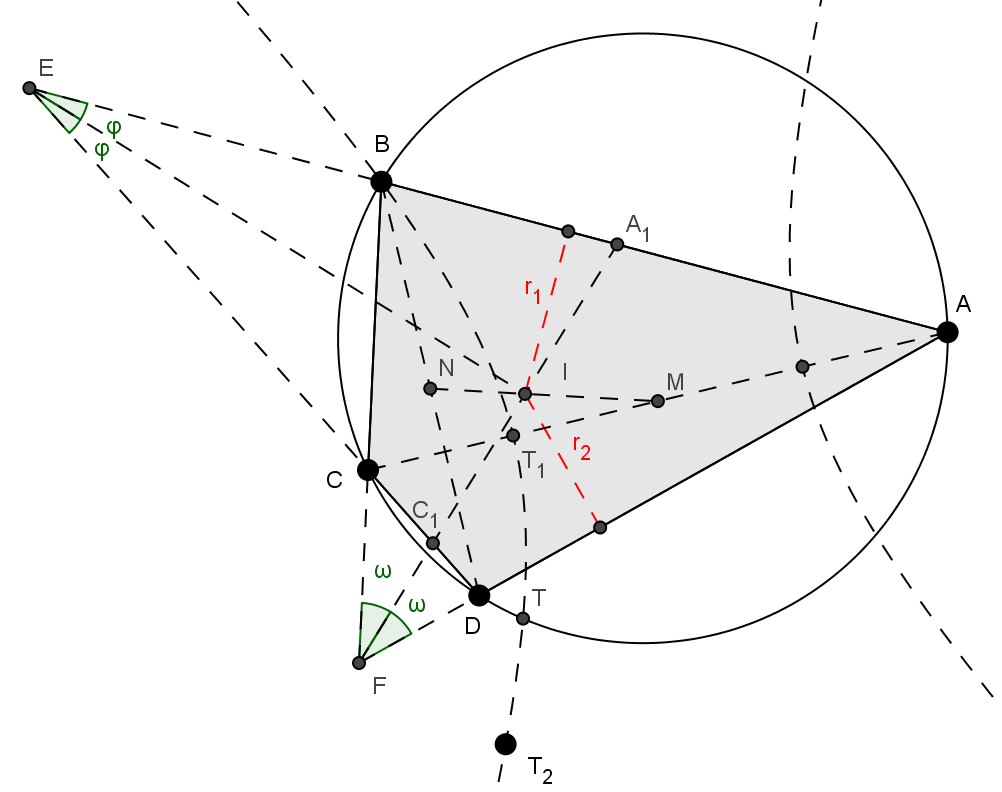}
\caption{Cyclic quadrilateral with $I_1\equiv I_2\equiv I$.
\label{newfig4}}
  \end{center}
  \vspace{-0.5cm}
\end{figure}

In the special case of a cyclic quadrilateral, the two incenters $I_1$ and $I_2$ coincide, but $r_1\neq r_2$. The following theorem has been proven in the distant past albeit in a different way:

\noindent
{\bf Theorem 5 (Based on Theorem \# 387 [3]).} In a cyclic quadrilateral $ABCD$, the incenters $I_1, I_2$ coincide with point $I$ on the segment MN of the Newton line (Figure~\ref{newfig4}). 

\noindent
{\bf Proof.} Since $I_1, I_2$ are located on the Newton line of $ABCD$, it is sufficient to show that the intersection $I$ of the two angle bisectors from $\angle E$ and $\angle F$ (Figure~\ref{newfig4}) lies on the Newton line $MN$ of $ABCD$. We find that $\angle ECF=\angle EIF +\phi +\omega$ and $\angle EAF=\angle EIF -\phi -\omega$ from which we get $\angle EIF = (\angle ECF + \angle EAF)/2 = 90^o$, so $\triangle EC_1A_1$ is isosceles and $I$ is the midpoint of $C_1A_1$.

Now let the diagonals be $AC=e$ and $BD=f$ and define $x=f/(e+f)$ and $y=e/(e+f)$, such that $x+y=1$. From the similar triangles $\triangle FAB\sim\triangle FCD$, $\triangle FAC\sim\triangle FBD$, and the angle bisector theorem, we find the proportions
\begin{equation}
\frac{AA_1}{A_1B} = \frac{FA}{FB} = \frac{e}{f} = \frac{FC}{FD} = \frac{CC_1}{C_1D} = \frac{y}{x} \, ,
\label{ratios1} 
\end{equation}
which imply that $A_1=xA+yB$, $C_1=xC+yD$, and finally for the midpoint $I$ of $C_1A_1$ that
\begin{equation}
I = (A_1+C_1)/2 = x(A+C)/2 + y(B+D)/2 = xM + yN \, ,
\label{ratios2} 
\end{equation}
which shows that $I$ lies on segment MN of the Newton line and divides it in a ratio of $MI:IN = y:x=e:f$, just as $A_1, C_1$ divide $AB, CD$, respectively.
\hfill\sq

\begin{figure}
\begin{center}
    \leavevmode
    \includegraphics[trim=.1cm .1cm .1cm .1cm, clip, angle=0,width=15 cm]{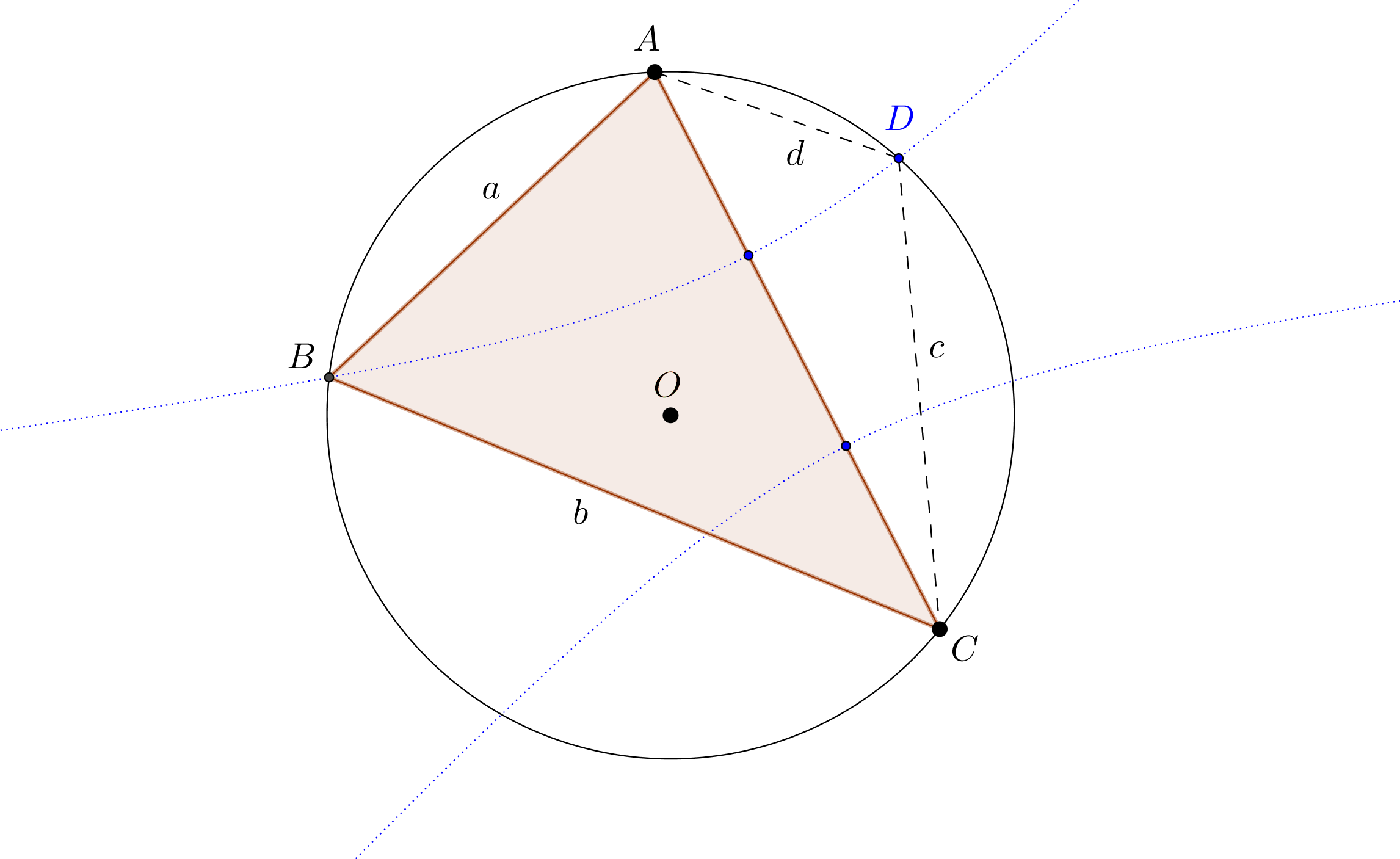}
\caption{Triangle $ABC$ inscribed in $\odot O$ and the bicentric quadrilateral $ABCD$ constructed by locating vertex $D$ on minor $\protect\arc{CA}$ as stated in Theorem 6. The hyperbola with foci $A$ and $C$ with one branch through vertex $B$ is shown by dotted lines. 
\label{fig3}}
  \end{center}
  \vspace{-0.5cm}
\end{figure}

\section{Bicentric Quadrilateral}

In the special case of a cyclic or tangential quadrilateral, we derive the conditions under which it is also tangential or cyclic, respectively, thus it is bicentric with a single incenter $I$ and inradius $r$. We prove the following two theorems:

\noindent
{\bf Theorem 6 (Cyclic Quadrilateral is Bicentric).} Consider $\triangle ABC$ inscribed in $\odot O$ (Figure~\ref{fig3}). Any point $D$ chosen on minor $\arc{CA}$ defines a cyclic quadrilateral $ABCD$ with side lengths $AB=a, BC=b, CD=c$, and $DA=d$. Of these quadrilaterals, there exists only one that is also tangential (therefore it is bicentric) to a single incircle $\odot I$. Its vertex $D$ lies at the point of intersection of minor $\arc{CA}$ with the branch of the hyperbola with foci $A$ and $C$ that passes through vertex $B$. An analogous property holds when point $D$ is chosen to lie on minor $\arc{AB}$ or on minor $\arc{BC}$.

\noindent
{\bf Proof.} Consider a hyperbola with foci $A$ and $C$ such that one of its branches passes through vertex $B$ and intersects minor $\arc{CA}$ at point $D$ (Figure~\ref{fig3}). Then by the geometric definition (the locus) of the hyperbola, we can write that
\begin{equation}
b-a = c-d \, ,
\label{hyp2} 
\end{equation}
where $a<b$ and $d<c$, as in Figure~\ref{fig3}. Similarly, in the case with $a>b$ and $d>c$, we can write that
\begin{equation}
a-b = d-c \, .
\label{hyp1} 
\end{equation}
Both equations are equivalent to 
\begin{equation}
a+c = b+d \, ,
\label{hyp5} 
\end{equation}
which implies that $ABCD$ is tangential (thus also bicentric).
\hfill\sq

\noindent
{\bf Theorem 7 (Tangential Quadrilateral is Bicentric).} Consider a tangential quadrilateral $ABCD$ with incenter $I$ (Figure~\ref{fig2}). Then $ABCD$ is cyclic (thus also bicentric) only if vertex $D$ lies on the hyperbola with foci $A$ and $C$ that passes through vertex $B$ (Figure~\ref{fig3}). An analogous property holds when any other vertex is chosen instead of $D$.

\noindent
{\bf Proof.} Consider the tangential quadrilateral $ABCD$ with incircle $\odot I$ shown in Figure~\ref{fig2}. By the tangency of its sides, equation~(\ref{in1}) is valid and it can be written in the form of equation~(\ref{hyp5}). We re-arrange terms in equation~(\ref{hyp5}) to obtain:
\begin{equation}
b-a = c-d \, .
\label{in12} 
\end{equation}
This equation defines a hyperbola with foci $A$ and $C$ one branch of which passes through vertices $B$ and $D$ (as shown in Figure~\ref{fig3} for the case $a < b, d<c$). If $D$ also lies on the circumcircle of $\triangle ABC$, then $ABCD$ is cyclic (thus also bicentric). 
\hfill\sq

\section{A Euclidean Construction of Vertex \boldmath{$D$}}

\begin{figure}
\begin{center}
    \leavevmode
    \includegraphics[trim=.1cm .1cm .1cm .1cm, clip, angle=0,width=15 cm]{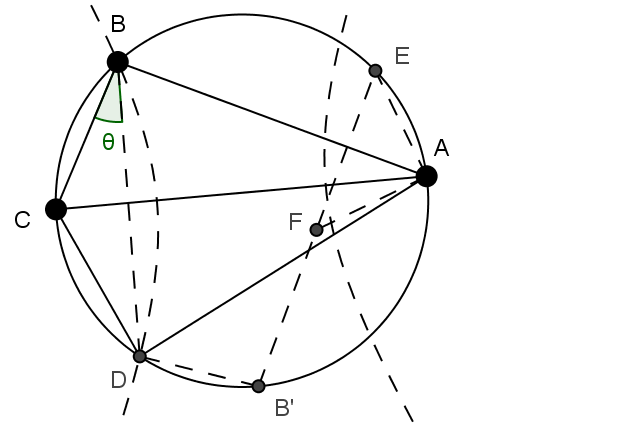}
\caption{A Euclidean construction of $D$ on the circumcircle of $\triangle ABC$ without using the hyperbola.
\label{fig7}}
  \end{center}
  \vspace{-0.5cm}
\end{figure}

Given $\triangle ABC$ inscribed in $\odot O(R)$ (as in Figure~\ref{fig7}), we construct on the lower $\arc{CA}$ of $\odot O(R)$ the point $D$, without using the hyperbola mentioned above, such that the quadrilateral $ABCD$ is convex and tangential with $AB+CD=BC+DA$.

\noindent
{\bf Construction.} Let $AB > BC$, as in Figure~\ref{fig7}. From the midpoint $B^\prime$ of the lower $\arc{CA}$, we draw a perpendicular to $AB$ that meets the circle at $E$, and a perpendicular to $EA$ at $A$ that meets $EB^\prime$ at $F$. On the minor $\arc{CB^\prime}$, we locate the required point $D$ such that $DB^\prime = FA$.

\noindent
{\bf Proof.} Since $\arc{BCB^\prime}$ + $\arc{AE}$ = $180^o$, then $2A + B$ + $\arc{AE} = A+B+C$ or $\arc{AE} = C - A$. Using this result, we find that 
\begin{equation}
 \begin{array}{ll} 
DB^\prime & = FA     \\
                  & = AE\tan\frac{B}{2} \\
                  & = 2R\sin\frac{C-A}{2}\tan\frac{B}{2} \, . \\
         \end{array} 
\label{con1}
\end{equation}

Next we define $\angle CBD = \theta$ and we find that
\begin{equation}
 \begin{array}{ll} 
      & AB +CD = BC + DA     \\
 \iff & 2R\sin C + 2R\sin\theta = 2R\sin A + 2R\sin(B-\theta)  \\
 \iff & 2R\sin\frac{C-A}{2}\cos\frac{C+A}{2}  = 2R\sin(\frac{B}{2}-\theta)\cos\frac{B}{2}\\
 \iff & 2R\sin(\frac{B}{2}-\theta) = 2R\sin\frac{C-A}{2}\tan\frac{B}{2}\\
 \iff & DB^\prime = 2R\sin\frac{C-A}{2}\tan\frac{B}{2}\, , \\
         \end{array} 
\label{con2}
\end{equation}
as was also found in equation~(\ref{con1}).\hfill\sq

\section{A Quadrilateral with \boldmath{$I_1\equiv I_2$}}

\noindent
{\bf Theorem 8.} If $I_1\equiv I_2$, then the quadrilateral is tangential, or cyclic, or bicentric.

\begin{figure}
\begin{center}
    \leavevmode
    \includegraphics[trim=.1cm .1cm .1cm .1cm, clip, angle=0,width=15 cm]{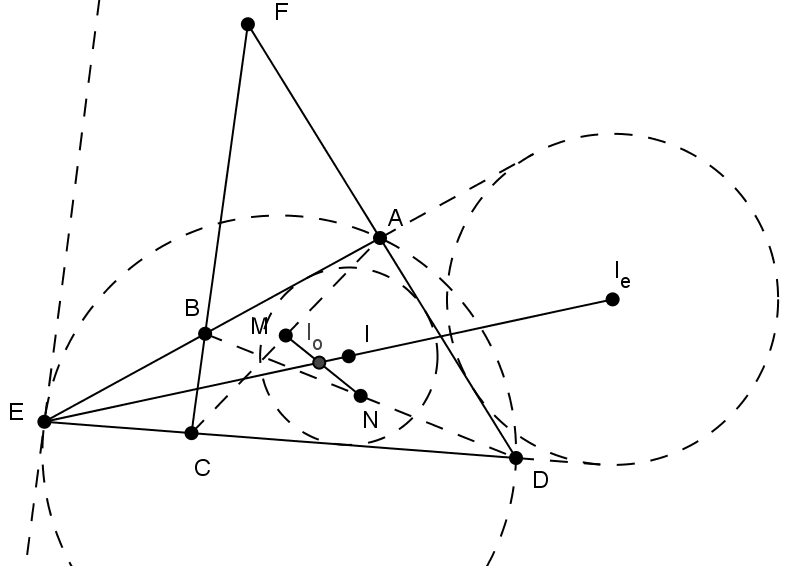}
\caption{Quadrilateral $ABCD$ with $I_1\equiv I_2\equiv I_o$. 
\label{fig8}}
  \end{center}
  \vspace{-0.5cm}
\end{figure}

\noindent
{\bf Proof.} Using barycentric coordinates in the basic $\triangle EDA$, let $E=(1:0:0)$, $D=(0:1:0)$, $A=(0:0:1)$, $DA=a$, $AE=b$, and $ED=c$. Also let $I_1\equiv I_2\equiv I_o$ on the bisector of $\angle DEA$ with barycentric coordinates $I_o=(k:b:c)$. Finally, let $I$ and $I_e$ be the incenter and $E$-excenter of $\triangle EDA$, respectively (Figure~\ref{fig8}).

Since $F=(0:m:1)$ and $B=(1:0:n)$, we obtain the equation $mnx+y -mz=0$ for the line $BC$ with point at infinity $P_\infty=(1+m:-m-mn:-1+mn)$, $C=(1:-mn:0)$, as well as the midpoints of the diagonals $M=(1:-mn:1-mn)$ and $N=(1:1+n:n)$.

Since the points $M, N,$ and $I_o$ are collinear, the determinant $det(M, N, I_o)$ must be zero from which we find that
\begin{equation}
n = \frac{b+c-k}{(b-c+k) + (b-c-k)m}\, .
\label{big1} 
\end{equation}

If a point $(x_o:y_o:z_o)$ is equidistant from the lines $a_ix+b_iy+c_iz=0$ ($i=1, 2$), the following equality holds [5]:
\begin{equation}
S_2 (a_1x_o+b_1y_o+c_1z_o)^2 = S_1 (a_2x_o+b_2y_o+c_2z_o)^2\, ,
\label{big2} 
\end{equation}
where
\begin{equation}
S_i = (b^2+c^2-a^2)(b_i - c_i)^2 + 
         (c^2+a^2-b^2)(c_i - a_i)^2 + 
         (a^2+b^2-c^2)(a_i - b_i)^2 \, ,
\label{big3} 
\end{equation}
and $i=1, 2$. Point $I_o$ is equidistant from the lines $BC$ (given above) and $DA$ ($x=0$).
Applying equation~(\ref{big2}) to $I_o$ and using equation~(\ref{big1}) to eliminate $n$, we find that
\begin{equation}
(1+m)(k^2-a^2)(b-cm)\left[b(b-c+k) - (b-c-k)cm\right] = 0\, .
\label{big4} 
\end{equation}
The solutions of this equation can be classified as follows:
\begin{itemize}
\item[(a)] $m=-1$ is rejected because then $F=(0:-1:1)$ becomes a point at infinity.
\item[(b)] $k=\pm a$ implies that $I_o$ coincides with $I$ or $I_e$, thus $\odot I_1(r_1)$ is tangent to a third side and by Theorem 4 this quadrilateral is tangential.
\item[(c)] $b=cm$ and $[b(b-c+k) - (b-c-k)cm]=0$ both imply that $P_\infty=(b^2-c^2:-b^2:c^2)$ which means that $BC$ is parallel to the tangent to the circumcircle of $\triangle EDA$ at point $E$, thus this quadrilateral is cyclic.
\end{itemize}
Solutions (b) and (c) of equation~(\ref{big4}) are the only ones that are valid. When one solution from (b) and one solution from (c) are simultaneously valid, then the quadrilateral is bicentric. This completes the proof.\hfill\sq





\vskip0.25in
\noindent
\author{Nikolaos Dergiades: I. Zanna 27, Thessaloniki 54643, Greece} \\
{\it Email address}: {\tt ndergiades@yahoo.gr}\\

\noindent
\author{Dimitris M. Christodoulou: Department of Mathematical Sciences, University of Massachusetts \\ Lowell, Lowell, MA 01854, USA} \\
{\it Email address}: {\tt dimitris\_christodoulou@uml.edu}

\end{document}